# A Strong Solution of the Navier-Stokes Equation on $L_p(\Omega)$

Maoting Tong

Department of Mathematical Science, Xi'anJiaotong-Liverpool University, Suzhou, 215123, P.R.China, E-mail address: Maoting.Tong18@student.xjtlu.edu.cn

Daorong Ton*

Hohai University, Nanjing, 210098, P.R.China, E-mail address: 1760724097@qq.com, Current address : 1-3306 Moonlight Square, Nanjing, 210036, P.R.China.

*- The corresponding author

## Abstract

In this paper we prove that the Navier-Stokes initial value problem (1.1)-(1.4) has a unique smooth local strong solution $u(t,x)$ and $p(t,x)$ if $f(t) \in DL_p(\Omega)$ is Hölder continuous about $t$ on $(0,T]$ and the initial value $u_0 \in D((-\Delta)^{\frac{1}{2}}) \cap (C^\infty(\Omega))^3$ satisfies $\frac{\partial u_{0i}}{\partial x_j} = 0 (i \neq j)$.

（1）**Introduction**

The Navier-Stokes initial value problem can be written in the following form as

$$\begin{cases} \frac{\partial u}{\partial t} = \Delta u - \nabla p - (u \bullet \nabla)u + f, x \in \Omega, t \in (t_0, \infty) & \text{...............(1.1)} \\ \nabla \bullet u \equiv div u = 0 & \text{...............(1.2)} \\ u\big|_{\partial\Omega} = 0, t \in (0, \infty) & \text{...............(1.3)} \\ u\big|_{t=t_0} = u_0, x \in \Omega & \text{...............(1.4)} \end{cases}$$

where $\Omega$ is a bounded domain in $R^3$ with smooth boundary $\partial\Omega$ of class $C^3$, $u =$

$u(t,x) = (u_1(t,x), u_2(t,x), u_3(t,x))$ is the velocity field, $u_0 = u_0(x)$ is the initial velocity, $p = p(t,x)$ is the pressure, $f = f(t,x)$ is the external force. In [17] we proved that the Navier-Stokes initial value problem (1.1)-(1.4) has a unique smooth local strong solution on $L_2(\Omega)$ if the initial velocity $u_0$ and the external force $f$ satisfy some conditions. In this paper we will extend this result to $L_p(\Omega)$.

Let $L_p(\Omega)$ ($2 \le p < \infty$) be the Banach space of real vector functions in $L^p(\Omega)$ with the inner product defined in the usual way. That is

$$L_p(\Omega) = \{u : \Omega \to R^3, u = (u_1, u_2, u_3), u_i \in L^p(\Omega)(i=1,2,3)\}.$$

For $u = (u_1, u_2, u_3) \in L_p(\Omega)$, we define the norm

$$\|u\|_{L_p(\Omega)} = (\sum_{i=1}^{3} \|u_i\|_{L^p(\Omega)}^p)^{\frac{1}{p}}.$$

The set of all real vector functions $u$ such that $div\ u = 0$ and $u \in C_0^\infty(\Omega)$ is denoted by $C_{0,\sigma}^\infty(\Omega)$. Let $DL_p(\Omega)$ be the closure of $C_{0,\sigma}^\infty(\Omega)$ in $L_p(\Omega)$. If $u \in C^\infty(\Omega)$ then $u \in DL_2(\Omega)$ implies $div\ u = 0$. (see p.270 in [4]). We have

$$C_{0,\sigma}^\infty(\Omega) \subset C_0^\infty(\Omega) \subset L_p(\Omega) \text{ and } DL_p(\Omega) = \overline{C_{0,\sigma}^\infty} \subset L_p(\Omega) \subset L_2(\Omega),$$

$$DL_p(\Omega) \subset L_p(\Omega) = W^{0,p}(\Omega), \quad \|\bullet\|_{DL_p(\Omega)} = \|\bullet\|_{L_p(\Omega)},$$

and

$$L_p(\Omega) = DL_p(\Omega) \oplus DL_p(\Omega)^\perp. \tag{2}$$

From [4] and [5] we see that $DL_p(\Omega)^\perp = \{\nabla h; h \in W^{1,p}(\Omega)\}$. Let $P$ be the orthogonal projection from $L_p(\Omega)$ onto $DL_p(\Omega)$. By applying $P$ to (1.1) and taking account of the other equations, we are let the following abstract initial value problem, $Pr.\text{II}$

$$\begin{cases} \dfrac{du}{dt} = P\Delta u + Fu + Pf, t \in (t_0, \infty) & \text{...............................................................(3.1)} \\ u|_{t=0} = u_0, x \in \Omega & \text{.............................................................................(3.2)} \end{cases}$$

where $Fu = -P(u \bullet \nabla)u$.

We consider equation (3.1) in integral form Pr.III

$$u(t) = e^{tP\Delta}u_0 + \int_{t_0}^{t} e^{(t-s)P\Delta}(Fu(s) + Pf)ds. \qquad (4)$$

For $u = (u_1, u_2, u_3) \in L_p(\Omega)$ we define $\Delta u = (\Delta u_1, \Delta u_2, \Delta u_3)$ and $\nabla u = (\dfrac{\partial u_1}{\partial x_1}, \dfrac{\partial u_2}{\partial x_2}, \dfrac{\partial u_3}{\partial x_3})$ where $\nabla$ is considered as an operator operating on vector functions. Since the operator $-\nabla = -\sum_{i=1}^{3} \dfrac{\partial^2}{\partial x_i^2}$ is strongly elliptic of order 2. From Theorem 7.3.6 in [12] $\Delta$ is the infinitesimal generator of an analytic semigroup of contractions on $L^p(\Omega)$ with $D(\Delta) = W^{2,p}(\Omega) \cap W_0^{1,p}(\Omega)$. Hence $\Delta$ is also the infinitesimal generator of an analytic semigroup of contraction on $L_p(\Omega)$ with $D(\Delta) = W_{2,p}(\Omega) \cap W_{1,p,0}(\Omega)$, where $W_{2,p}(\Omega)$ and $W_{1,p,0}(\Omega)$ are the Sobolev spaces of vector value in $W^{2,p}(\Omega)$ and $W_0^{1,P}(\Omega)$ respectively. We will prove that $\Delta$ is also the infinitesimal generator of an analytic semigroup of contraction on $DL_p(\Omega)$. This is a key point for our purpose.

**(2) Some lemmas**

**Lemma 1.** If $u \in C^\infty(\Omega)$ then $u \in DL_p(\Omega)$ implies $div\ u = 0$.

Proof. Suppose that $u = (u_1, u_2, u_3) \in C^\infty(\Omega) \cap DL_p(\Omega)$. Then there exists a sequence $\{u^n \in C_{0,\sigma}^\infty(\Omega) : n = 1, 2, ...\}$ such that $u = \lim_{n \to \infty} u^n$, that is $\lim_{n \to \infty} \|u - u^n\|_{DL_p(R^3)} = 0$ uniformly on $\Omega$ and $\lim_{n \to \infty} \|u_i - u_i^n\|_{L^p(\Omega)} = 0$ uniformly on $\Omega$ for $i = 1, 2, 3$. And so $\lim_{n \to \infty} |u_i(x) - u_i^n(x)| \overset{a.e.}{=} 0$ and $\lim_{n \to \infty} u_i^n(x) \overset{a.e.}{=} u_i(x)$ uniformly on $\Omega$ for $i = 1, 2, 3$. (see Theorem 1.39 in [15]）It from the proof of Theorem 7.11 in [14] follows that

$$divu = \sum_{i=1}^{3}\frac{\partial u_i}{\partial x_i} = \sum_{i=1}^{3}\frac{\partial(\lim_{n\to\infty}u_i^n)}{\partial x_i} = \sum_{i=1}^{3}\left[\lim_{x_i\to x_{i0}}\frac{\lim_{n\to\infty}u_i^n(x_i)-\lim_{n\to\infty}u_i^n(x_{i0})}{x_i-x_{i0}}\right]$$

$$=\sum_{i=1}^{3}\left[\lim_{n\to\infty}\lim_{x_i\to x_{io}}\frac{u_i^n(x_i)-u_i^n(x_{i0})}{x_i-x_{i0}}\right] = \sum_{i=1}^{3}\left[\lim_{n\to\infty}\frac{\partial u_i^n}{\partial x_i}\right] = \lim_{n\to\infty}\sum_{i=1}^{3}\left[\frac{\partial u_i^n}{\partial x_i}\right] = 0.$$

This is to say $div\ u = 0$. □

**Lemma 2.** For every $u \in L_p(\Omega)$, $div\ u = 0$ if and only if $div\ (\lambda I - \Delta)u = 0$ for

$$\lambda \in \Sigma_\vartheta = \{\lambda : \vartheta - \pi < \arg\lambda < \pi - \vartheta, |\lambda| \geq r\}$$

where $0 < \vartheta < \frac{1}{2}$.

Proof. It is similar to the proof of lemma 1 in [17].

**Lemma 3.** (1.5.12 in[8]) Let $\{T(t): t \geq 0\}$ be a $C_0$-semigroup on a Banach space $X$. If $Y$ is a closed subspace of $X$ such that $T(t)Y \subset Y$ for all $t \geq 0$, i.e., if $Y$ is $T(t)_{t\geq 0}$-invariant, then the restrictions

$$T(t)_| := T(t)_{|Y}$$

form a $C_0$-semigroup $\{T(t)_| : t \geq 0\}$, called the subspace semigroup on the Banach space $Y$.

**Lemma 4.** (Proposition 2.2.3 in [8]) Let $(A, D(A))$ be the generator of a $C_0$-semigroup $\{T(t): t \geq 0\}$ on a Banach space $X$ and assume that the restricted semigroup (subspace semigroup) $\{T(t)_| : t \geq 0\}$ is a $C_0$-semigroup on some $(T(t))_{t\geq 0}$ – invariant Banach space $Y \to X$. Then the generator of $\{T(t)_| : t \geq 0\}$ is the part $(A_|, D(A_|))$ of $A$ in $Y$.

**Lemma 5.** The operator $\Delta_{|DL_p(\Omega)}$ is the infinitesimal generator of an analytic semigroup of contractions on $DL_p(\Omega)$.

Proof. We have know that $\Delta$ is a infinitesimal generator of an analytic semigoup of contraction on $L_p(\Omega)$. Let $\{T(t)|t \geq 0\}$ be the restriction of the analytic semigroup generated by $\Delta$ on $L_p(\Omega)$ to the real axis. Then $\{T(t)|t \geq 0\}$ is a $C_0$ semigroup of contractions. We have already noted that $DL_p(\Omega)$ is a closed subspace of $L_p(\Omega)$. We want to show that $DL_p(\Omega)$ is $T(t)_{t \geq 0}$ – invariant.

For every $u \in L_p(\Omega)$ with $div\ u = 0$ and $\lambda \in \rho(\Delta) \cap \Sigma_\vartheta = \{\lambda : \vartheta - \pi < \arg \lambda < \pi - \vartheta, |\lambda| \geq r\}$ we have $(\lambda I - \Delta)[R(\lambda : \Delta)u] = u$ where $\Sigma_\vartheta$ is the same as in the proof of lemma 2. From Lemma 2 it follows that $div\ R(\lambda : \Delta)u = 0$. That is to say that $R(\lambda : \Delta)$ is preserving divergence-free for $\lambda \in \rho(\Delta) \cap \Sigma_\vartheta$. From Theorem 2.5.2 (c) in [12] it follows that $\rho(\Delta) \supset R^+$, and so $\rho(\Delta) \cap \Sigma_\vartheta \supset \{\lambda : \lambda \geq r\}$. Hence $R(\lambda : \Delta)$ is preserving divergence-free for every $\lambda \geq r$. Let $u \in DL_p(\Omega)$ then there exists a sequence $u_n$ such that $\lim_{n \to \infty} u_n = u$ and $div\ u_n = 0$ for $n = 1, 2, \ldots$. Since $R(\lambda : \Delta)$ is bounded and so is continuous. Hence $\lim_{n \to \infty} R(\lambda : \Delta)u_n = R(\lambda : \Delta)u$ and $div\ R(\lambda : \Delta)u_n = 0$ for every $\lambda \geq r$. Therefore $R(\lambda : \Delta)u \in DL_p(\Omega)$ for every $\lambda \geq r$. It follows that $DL_p(\Omega)$ is $R(\lambda : \Delta)$ -invariant for every $\lambda \geq r$. Now Theorem 4.5.1 in [12] implies that $DL_p(\Omega)$ is $T(t)_{t \geq 0}$ – invariant. From Lemma 3 and Lemma 4 it follows that $\Delta_{|DL_p(\Omega)}$ is the infinitesimal generator of the $C_0$ semigroup $\{T(t)_{|DL_p(\Omega)} : t \geq 0\}$ of contractions on $DL_p(\Omega)$.

By the proof of Theorem 2.5.2 in [12] we see that the analytic semigroup generated by infinitesimal generator $\Delta$ on $L_p(\Omega)$ has a representation

$$T(z) = T(t) + \sum_{n=1}^{\infty} \frac{T^{(n)}(t)}{n!}(z-t)^n$$

for $z = t + is \in \Delta_1 = \{z : |\arg z| \langle \arctan 1/Ce\}$ where $C$ is a constant with $\|T'(t)\| = \|\Delta T(t)\| \leq C/t$ by Theorem 2 5.2(d) in [12]. For every $u \in DL_p(\Omega)$ we have

$$T(z)u$$
$$= T(t)u + \sum_{n=1}^{\infty} \frac{T^{(n)}(t)}{n!}(z-t)^n u$$
$$= T(t)u + \sum_{n=1}^{\infty} \frac{\left[\Delta_{|DL_2(\Omega)}T(\tfrac{1}{n})\right]^n}{n!}(is)^n u$$
$$= T(t)u + \lim_{k\to\infty} \sum_{n=1}^{k} \frac{\left[\Delta_{|DL_2(\Omega)}T(\tfrac{1}{n})\right]^n}{n!}(is)^n u$$
$$= T(t)u - \lim_{k\to\infty} \sum_{n=1}^{k} s^{2n} \frac{\left[\Delta_{|DL_2(\Omega)}T(\tfrac{1}{n})\right]^{2n}}{(2n)!} u .$$

We used the lemma 2.4.5 in [12] in the second step, and the fact that subsequence converges to the same limit in the fourth step. Since $DL_p(\Omega)$ is $\Delta_{|DL_p(\Omega)}$ – invariant and $T(t)_{t\geq 0}$ – invariant, hence

$$\left[\Delta_{|DL_p(\Omega)}T(\tfrac{1}{n})\right]^{2n} u \in DL_p(\Omega)$$

and so

$$\sum_{n=1}^{k} s^{2n} \frac{\left[\Delta_{|DL_p(\Omega)}T(\tfrac{1}{n})\right]^{2n}}{(2n)!} u \in DL_p(\Omega).$$

Finally we have

$$\lim_{k\to\infty} \sum_{n=1}^{k} s^{2n} \frac{\left[\Delta_{|DL_p(\Omega)}T(\tfrac{1}{n})\right]^{2n}}{(2n)!} u \in DL_p(\Omega)$$

because $DL_p(\Omega)$ is closed. Therefore $T(z)u \in DL_p(\Omega)$ for $z \in \Delta_1$. This is to say that $DL_p(\Omega)$ is $T(z)_{z\in\Delta_1}$ – invariant and $\Delta_{|DL_p(\Omega)}$ is a infinitesimal generator of the analytic semigroup $\{T(z): z \in \Delta_1\}$ on $DL_p(\Omega)$. □

Suppose that $-A$ is the infinitesimal generator of an analytic semigroup $T(t)$ on a Banach space $X$. From the results of section 2.6 in [12] we can define the fraction powers $A^\alpha$ for $0 \leq \alpha \leq 1$ and $A^\alpha$ is a closed linear invertible operator with domain $D(A^\alpha)$ dense in $X$ and $A^{-\alpha}$ is bounded. $D(A^\alpha)$ equipped with the norm $\|x\|_\alpha = \|A^\alpha x\|$ is a Banach space denoted by $X_\alpha$. It is clear that $0 < \alpha < \beta$ implies $X_\alpha \supset X_\beta$ and that the embedding of $X_\beta$ into $X_\alpha$ is continuous. If $-A = \Delta$ and $\gamma > 3/4 > 1/2$ then $X_\gamma \subset X_{\frac{1}{2}}$ and $D((-\Delta)^\gamma) \subset D((-\Delta)^{1/2})$, so the norms $\|\bullet\|_{DL_p(R^3)_\gamma}$ and $\|\bullet\|_{DL_p(R^3)_{\frac{1}{2}}}$ are equivalent (see p291 in [9]), i.e. there exists $L_0 > 0$ such that for any $u \in D((-\Delta)^\gamma)$

$$\|(-\Delta)^\gamma u\|_{DL_p(\Omega)} = \|u\|_{DL_p(\Omega)_\gamma} \le L_0 \|u\|_{DL_p(\Omega)_{1/2}}. \tag{5}$$

For $u \in D(\Delta)$ we have

$$\|\nabla u\|_{DL_p(\Omega)} = \|-\nabla u\|_{DL_p(\Omega)} = \|(-\Delta)^{1/2} u\|_{DL_p(\Omega)} = \|u\|_{DL_p(\Omega)_{1/2}}. \tag{6}$$

In [5] Giga proved the following result:

**Lemma 6.** (Lemma 2.2 in [5]) Let $0 \le \delta < 1/2 - n(1-r^{-1})/2$. Then

$$\|A^{-\delta} P(u, \nabla)v\|_{0,r} \le M \|A^\vartheta u\|_{0,r} \|A^\omega v\|_{0,r}$$

with some constant $M = M(\delta, \vartheta, \omega, r)$, provided $\delta + \vartheta + \omega \ge n/2r + 1/2$, $\vartheta > 0$, $\omega > 0$, $\omega + \delta > 1/2$.

From the Lemma 6 and the formula (5) we see that if take $n = 3, r = 2$ $\delta = 0$, $\vartheta = 3/4$ and $\omega = 3/4$, then

$$\|(u \bullet \nabla)v\|_{DL_p(\Omega)} = \|(u \bullet \nabla)v\|_{L_p(\Omega)} = \|(u \bullet \nabla)v\|_{W^{0,p}(\Omega)}$$
$$\le M \|(-\Delta)^\vartheta u\|_{DL_p(\Omega)} \|(-\Delta)^\omega v\|_{DL_p(\Omega)}$$
$$\le M L_0^2 \|u\|_{DL_p(\Omega)_{1/2}} \|v\|_{DL_p(\Omega)_{1/2}}$$

with some constant $M$ for any $u, v \in DL_p(\Omega)$. Hence we have

**Lemma 7.** Suppose that $u$, $v \in DL_p(\Omega)$ are velocity fields and $(u \bullet \nabla)v \in DL_p(\Omega)$, then

$$\|(u \bullet \nabla)v\|_{DL_p(\Omega)} \le M L_0^2 \|u\|_{DL_p(\Omega)_{1/2}} \|v\|_{DL_p(\Omega)_{1/2}}.$$

Now we introduce the main lemmas of this paper.

**Assumption (F).** Let $X = DL_p(\Omega)$ and $U$ be an open subset in $R^+ \times X_\alpha (0 \leq \alpha \leq 1)$

The function $f: U \to X$ satisfies the assumption (F) if for every $(t,u) \in U$ there is a neighborhood $V \subset U$ and constants $L \geq 0$, $0 < \vartheta \leq 1$ such that for all $(t_i, u_i) \in V (i = 1,2)$

$$\|f(t_1, u_1) - f(t_2, u_2)\|_X \leq L(|t_1 - t_2|^\vartheta + \|u_1 - u_2\|_\alpha). \tag{7}$$

**Lemma 8.** ( Theorem 6.3.1 in [12]) Let $-A$ is the infinitesimal generator of an analytic semigroup $T(t)$ on the Banach space $X = DL_p(\Omega)$ satisfying $\|T(t)\| \leq 1$ and assume further that $0 \in \rho(-A)$. If, $0 < \alpha < 1$ and $f$ satisfies the assumption $(F)$, then for every initial data $(t_0, u_0) \in U$ the initial value problem

$$\begin{cases} \dfrac{du(t)}{dt} + Au(t) = f(t, u(t)), t \in (t_0, \infty) & \text{...............(8.1)} \\ u(t_0) = u_0 & \text{...............(8.2)} \end{cases}$$

has a unique local solution

$$u \in C([t_0, T): DL_p(\Omega)) \cap C^1((t_0, T): DL_p(\Omega))$$

where $T = T(u_0)$.

In what follows we will need Banach lattice (see [3] ). A real vector space $G$ which is ordered by some order relation $\leq$ is called a vector lattice (or Riesz space) if any two elements $f, g \in G$ have a least upper bound, denoted by $f \vee g$, and a greatest lower bound, denoted by $f \wedge g$, and the following properties are satisfied:

(i) If $f \leq g$, then $f + h \leq g + h$ for all $f, g, h \in G,$

(ii) If $0 \leq f$, then $0 \leq tf$ for all $f \in G$ and $0 \leq t \in R$.

A Banach lattice is a real Banach space $G$ endowed with an ordering $\leq$ such that $(G, \leq)$ is a vector lattice and the norm is a lattice norm, that is $|f| \leq |g|$ implies $\|f\| \leq \|g\|$ for $f, g \in G$, where $|f| = f \vee (-f)$ is the absolute value of $f$ and $\|\bullet\|$ is the norm in $G$. In a Banach lattice $G$ we define for $f \in G$

$$f^+ := f \vee 0, \quad f^- := (-f) \vee 0.$$

The absolute value of $f$ is $|f| = f^+ + f^-$ and $f = f^+ - f^-$.

$$0 \le f \le g \Rightarrow |f| \le |g| \Rightarrow \|f\| \le \|g\|.$$

In what follows one will need the above formula.

In $L^p(\Omega)$ we define the order for $f, g \in L^p(\Omega)$

$$f \le g \Leftrightarrow f(x) \le g(x) \text{ for } a.e. \ x \in \Omega$$

and $(f \vee g)(x) := \max\{f(x), g(x)\}, \quad (f \wedge g)(x) := \min\{f(x), g(x)\}$

for $a.e.\ x \in \Omega$. Then $W^{1,p}(\Omega), L^p(\Omega), L_p(\Omega)$ and $DL_p(\Omega)$ are all Banach lattices. (see [1] p.148)

**Lemma 9.** Suppose that $u, v \in DL_p(\Omega)$ are divergence free satisfying $\dfrac{\partial u_i}{\partial x_j} = 0, \dfrac{\partial v_i}{\partial x_j} = 0$ $(i \ne j)$. Then $(u \bullet \nabla)v, \ (v \bullet \nabla)u \in DL_p(\Omega)$.

Proof. If $u, v \in DL_p(\Omega)$ are divergence free satisfying $\dfrac{\partial u_i}{\partial x_j} = 0, \dfrac{\partial v_i}{\partial x_j} = 0 (i \ne j)$. From 【5】 and 【6】 we have

$$\int_\Omega u \bullet \nabla h \, dx = 0, \int_\Omega v \bullet \nabla h \, dx = 0 \text{ for all } h \in W^{1,p}(\Omega).$$

That is $\int_\Omega \sum_{i=1}^3 u_i \dfrac{\partial h}{\partial x_i} dx = 0$. Since $L_p(\Omega)$ is a Banach lattice and $\nabla h \in DL_p(\Omega)^\perp$ for $h \in W^{1,p}(\Omega)$, then $u = u^+ - u^-, u^+, u^- \in DL_p(\Omega), \nabla h = (\nabla h)^+ - (\nabla h)^-$. From proposition 10.8 in [3] the lattice operations $\wedge$ and $\vee$ are continuous, we have

$$(\nabla h)^+ = \begin{pmatrix} \dfrac{\partial h}{\partial x_1} \\ \dfrac{\partial h}{\partial x_2} \\ \dfrac{\partial h}{\partial x_3} \end{pmatrix}^+ = \begin{pmatrix} (\dfrac{\partial h}{\partial x_1})^+ \\ (\dfrac{\partial h}{\partial x_2})^+ \\ (\dfrac{\partial h}{\partial x_3})^+ \end{pmatrix} = \begin{pmatrix} (\lim_{x_1 \to x_{10}} \dfrac{h(x_1,x_2,x_3)-h(x_{10},x_2,x_3)}{x_1-x_{10}})^+ \\ (\lim_{x_2 \to x_{20}} \dfrac{h(x_1,x_2,x_3)-h(x_1,x_{20},x_3)}{x_2-x_{20}})^+ \\ (\lim_{x_3 \to x_{30}} \dfrac{h(x_1,x_2,x_3)-h(x_1,x_2,x_{30})}{x_3-x_{30}})^+ \end{pmatrix}$$

$$= \begin{pmatrix} \lim_{x_1 \to x_{10}} \dfrac{h^+(x_1,x_2,x_3)-h^+(x_{10},x_2,x_3)}{x_1-x_{10}} \\ \lim_{x_2 \to x_{20}} \dfrac{h^+(x_1,x_2,x_3)-h^+(x_1,x_{20},x_3)}{x_2-x_{20}} \\ \lim_{x_3 \to x_{30}} \dfrac{h^+(x_1,x_2,x_3)-h^+(x_1,x_2,x_{30})}{x_3-x_{30}} \end{pmatrix} = \begin{pmatrix} \dfrac{\partial h^+}{\partial x_1} \\ \dfrac{\partial h^+}{\partial x_2} \\ \dfrac{\partial h^+}{\partial x_3} \end{pmatrix} = \nabla(h^+),$$

Similarly, $(\nabla h)^- = \nabla(h^-)$. $W^{1,p}(\Omega) \subset W^{0,p}(\Omega) = L_p(\Omega)$ are all Banach spaces. $h \in W^{1,p}(\Omega)$ implies $h^+, h^- \in W^{1,p}(R^3)$, $\int_\Omega u^+ \bullet \nabla(h^+) dx = 0$. Since $v$ is divergence free, $\sum_{i=1}^3 \dfrac{\partial v}{\partial x_i} = 0$, so $\dfrac{\partial v}{\partial x_i}$ are all bounded, $\left|\dfrac{\partial v}{\partial x_i}\right| \le L$ ($i=1,2,3$) for some constant $L > 0$. We have

$$0 = -L\int_\Omega \sum_{i=1}^3 u_i^+ (\dfrac{\partial h}{\partial x_i})^+ dx \le \int_\Omega (u^+ \bullet \nabla)v \bullet (\nabla h)^+ dx$$

$$= \int_\Omega \sum_{i=1}^3 u_i^+ \dfrac{\partial v_i}{\partial x_i} (\dfrac{\partial h}{\partial x_i})^+ dx \le L \int_\Omega \sum_{i=1}^3 u_i^+ \dfrac{\partial (h^+)}{\partial x_i} dx = 0$$

Hence $\int_\Omega (u^+ \bullet \nabla)v \bullet (\nabla h)^+ dx = 0$. Similarly, $\int_\Omega (u^+ \bullet \nabla)v \bullet (\nabla h)^- dx = 0$. So we have
$$\int_\Omega (u^+ \bullet \nabla)v \bullet \nabla h dx = \int_\Omega (u \bullet \nabla)v \bullet [(\nabla h)^+ - (\nabla h)^-] dx = 0.$$
Similarly, $\int_\Omega (u^- \bullet \nabla)v \bullet \nabla h dx = 0$. Therefore

$$\int_\Omega (u \bullet \nabla)v \bullet \nabla h dx = \int_\Omega ((u^+ - u^-) \bullet \nabla)v \bullet \nabla h dx = \int_\Omega (u^+ \bullet \nabla)v \bullet \nabla h dx - \int_\Omega (u^- \bullet \nabla)v \bullet \nabla h dx = 0$$

and so $(u \bullet \nabla)v \in DL_p(\Omega)$. Similarly $(v \bullet \nabla)u \in DL_p(\Omega)$. □

**(3) Main result**

Now we study the Navier-Stokes initial value problem (1.1)-(1.4).

A function $u$ which is differentiable almost everywhere on $[0,T]$ such that $u' \in L^1[0,T:DL_2(\Omega)]$ is called a strong solution of the initial value problem (1.1)=(1.4) if $u(0) = u_0$ and $u$ satisfies (1.1)-(1.4) a.e. on $[0,T]$. A function $u:[0,T] \to DL_2(\Omega)$ is a classical solution of (3.1)-(3.2) on $[0,T]$ if $u$ is continuous on $[0,T]$, continuous differentiable on $(0,T)$, $u \in D(\Delta)$ for $t \in (0,T)$ and (3.1)-(3.2) are satisfied on $[0,T]$.

Let $H([0,T]; DL_p(\Omega)_{1/2})$ denote the space of all Hölder continuous functions $u(t,x)$ on $[0,T]$ with different exponents in $(0,1]$ and with divergence free functions values satisfying $\frac{\partial u_i}{\partial x_j} = 0$ $(i \neq j)$ in the Banach space $DL_p(\Omega)_{\frac{1}{2}}$. Then from lemma 9 for any $u \in H([0,T]; DL_p(\Omega)_{1/2})$ and any $t_1, t_2 \in [0,T]$, $(u(t_1) \bullet \nabla)u(t_2) \in DL_p(\Omega)$; and for any $u_1, u_2 \in H([0,T]; DL_p(\Omega)_{1/2})$ and any $t \in [0,T]$, $(u_1(t) \bullet \nabla)u_2(t) \in DL_p(\Omega)$. In the following we will use these facts. Let

$$H = \{u(t) : u \in H([0,T]; DL_p(\Omega)_{1/2}), t \in (0,T]\}$$

$H$ is a subset of $DL_p(\Omega)_{1/2}$ which consists of function values of all functions in $H([0,T]; DL_p(\Omega)_{1/2})$. Let $u_k(t,x) = (k_1, k_2, k_3)$ $(t \in [0,T], x \in \Omega, k_i \in R, i = 1,2,3))$. Then $u_k \in H([0,T]; DL_p(\Omega)_{1/2})$ and $u_k(t) \in H$ for all $k_1, k_2, k_3 \in R$ and all $t \in [0,T]$. Suppose that $u_i(x_i) \in DL_p(\Omega)(i = 1,2,3)$ satisfying $\sum_{i=1}^{3} \frac{\partial u_i(x_i)}{\partial x_i} = 0$. Let $u(x) = (u_1(x_1), u_2(x_2), u_3(x_3))$ and $u(t,x) \equiv u(x)$ for $t \in [0,T]$. Then $u(t,x) \in H([0,T]; DL_p(\Omega)_{1/2})$ and $u(x) \in H$. Hence $H$ is not empty. Take the open kernel $H^0$ of $H$ in $DL_p(\Omega)_{1/2}$. $H^0$ is also not empty. It is clear that $H^0$ is a open subset in $DL_p(\Omega)_{1/2}$. The bilinear form $(v \bullet \nabla)u$ on $H^0$ takes value in

$DL_p(\Omega)$. Let

$$U = (0,T) \times H^0.$$

Then $U$ is an open subset of $[0,T] \times DL_p(\Omega)_{1/2}$. If $u_0 \in H^0$, that is, there exist $u \in H([0,T]; DL_p(\Omega)_{1/2})$ and $t \in [0,T]$ such that $u_0 = u(t)$. Let $u_0(t) \equiv u_0$ for all $t \in [0,T]$. It is easy to see that $u_0(t) \in H([0,T]; DL_p(\Omega)_{1/2})$. Hence $u_0$ is also a value of another function. That is to say that a function $u(t) \in H^0$ can be value of different functions in $H([0,T]; DL_p(\Omega)_{1/2})$.

**Theorem**. The Navier-Stokes initial value problem (1.1)-(1.4) has a unique smooth local strong solution $u(t,x)$ and $p(t,x)$ if the following condition are satisfied

(1) $f \in DL_p(\Omega)$ is Hölder continuous about $t$ on $(0,T]$,

(2) The initial value $u_0 \in D((-\Delta)^{\frac{1}{2}}) \cap (C^\infty(\Omega))^3$ satisfies $\dfrac{\partial u_{0i}}{\partial x_j} = 0 (i \neq j)$.

Proof. ( Step1 ) First, $F(t,u(t)) = -(u(t) \bullet \nabla)u(t)$ is a function : $U \to DL_p(\Omega)$ because $(u(t) \bullet \nabla)u(t) \in DL_p(\Omega)$ according to the definition of $U$.

We will find that by incorporating the divergence-free condition we can remove the pressure term from our equation. (see p. $271^3$ in 【4】, p. $234_6$ and p. $239_9$ in 【13】) In fact, from $DL_p(\Omega)^\perp = \{\nabla h; h \in W^{1,p}(\Omega)\}$ we see that $\nabla p \in DL_p(\Omega)^\perp$ and so $P\nabla p = 0$. For $u \in DL_p(\Omega)$ we have $\Delta u \in DL_p(\Omega)$ because to Lemma 5. Hence by applying $P$ to the equation (1.1) we have $P\Delta u = \Delta u$. It follows from $(u \bullet \nabla)u \in DL_p(\Omega)$ that $P(u \bullet \nabla)u = (u \bullet \nabla)u$. Therefore we can first rewrite (1) into an abstract initial value problem on $DL_p(\Omega)$

$$\begin{cases} \dfrac{du}{dt} = \Delta u + F(t,u(t)), t \in [t_0,T] \\ u\big|_{t=t_0} = u_0, x \in \Omega \end{cases} \qquad (9)$$

where $F(t,u(t)) = -(u \bullet \nabla)u$. From Lemma 4 $\Delta_{|DL_p(\Omega)}$ is the generator of an analytic semigroup $T(t)$ of contraction on $DL_p(\Omega)$. So $\|T(t)\| \leq 1$. From Theorem 2.5.2(c) in 【15】 $0 \in \rho(\Delta)$.

(Step2) If $u(t)$ is Hölder continuous about $t$ on $[0,T]$ in $DL_p(\Omega)_{1/2}$, so there is a constant $C$ and $0 < \beta \leq 1$ such that

$$\|u(t_1,x) - u(t_2,x)\|_{DL_p(R^3)_{1/2}} \leq C|t_1 - t_2|^\beta \text{ for } t_1, t_2 \in [0,T]. \tag{10}$$

For any $(t_1, u_1(t_1)), (t_2, u_2(t_2)) \in U$ we have $(u_1(t_1) \bullet \nabla)u_1(t_1), (u_2(t_2) \bullet \nabla)u_2(t_2) \in DL_p(\Omega)$

$$\begin{aligned}
&\|(u_1(t_1) \bullet \nabla)u_1(t_1) - (u_2(t_1) \bullet \nabla)u_2(t_1)\|_{DL_p(\Omega)} \\
&= \|(u_1(t_1) \bullet \nabla)u_1(t_1) - (u_1(t_1) \bullet \nabla)u_2(t_1)\|_{DL_p(\Omega)} \\
&+ \|(u_1(t_1) \bullet \nabla)u_2(t_1) - (u_2(t_1) \bullet \nabla)u_2(t_1)\|_{DL_p(\Omega)} \\
&\leq \|(u_1(t_1) \bullet \nabla)(u_1(t_1) - u_2(t_1))\|_{DL_p(\Omega)} \\
&+ \|[(u_1(t_1) - u_2(t_1)) \bullet \nabla]u_2(t_1)\|_{DL_p(\Omega)} \\
&\leq ML_0^2 (\|u_1(t_1)\|_{DL_2(\Omega)_{1/2}} \|(u_1(t_1) - u_2(t_1))\|_{DL_p(\Omega)_{1/2}} \\
&+ \|u_1(t_1) - u_2(t_1)\|_{DL_p(\Omega)_{1/2}} \|u_2(t_2)\|_{DL_p(\Omega)_{1/2}}) \\
&= ML_0^2 \left( \|u_1(t_1)\|_{DL_p(\Omega)_{1/2}} + \|u_2(t_2)\|_{DL_p(\Omega)_{1/2}} \right) \|u_1(t_1) - u_2(t_2)\|_{DL_p(\Omega)_{1/2}} \\
&\leq ML_0^2 \left( \|u_1(t_1)\|_{DL_p(\Omega)_{1/2}} + \|u_2(t_2)\|_{DL_p(\Omega)_{1/2}} \right) \\
&\left( \|u_1(t_1) - u_2(t_2)\|_{DL_p(\Omega)_{1/2}} + \|u_2(t_2) - u_2(t_1)\|_{DL_p(\Omega)_{1/2}} \right) \\
&\leq ML_0^2 \left( \|u_1(t_1)\|_{DL_p(\Omega)_{1/2}} + \|u_2(t_2)\|_{DL_p(\Omega)_{1/2}} \right) \\
&\left( \|u_1(t_1) - u_2(t_2)\|_{DL_p(\Omega)_{1/2}} + C_1|t_1 - t_2|^{\beta_1} \right).
\end{aligned} \tag{11}$$

We used lemma 6 in the above third step. For any $u(t_1), u(t_2) \in DL_p(\Omega)_{1/2}$ we have

$$\|(u(t_1)\bullet\nabla)u(t_1)-(u(t_2)\bullet\nabla)u(t_2)\|_{DL_p(\Omega)}$$
$$\leq \|(u(t_1)\bullet\nabla)u(t_1)-(u(t_1)\bullet\nabla)u(t_2)\|_{DL_p(\Omega)} + \|(u(t_1)\bullet\nabla)u(t_2)-(u(t_2)\bullet\nabla)u(t_2)\|_{DL_p(\Omega)}$$
$$= \|(u(t_1)\bullet\nabla)(u(t_1)-u(t_2))\|_{DL_p(\Omega)} + \|[(u(t_1)-u(t_2))\bullet\nabla]u(t_2)\|_{DL_p(\Omega)}$$
$$\leq ML_0^2(\|u(t_1)\|_{DL_p(\Omega)_{1/2}}\|u(t_1)-u(t_2)\|_{DL_p(\Omega)_{1/2}} + \|u(t_1)-u(t_2)\|_{DL_p(\Omega)_{1/2}}\|u(t_2)\|_{DL_p(\Omega)_{1/2}})$$
$$= ML_0^2(\|u(t_1)\|_{DL_p(\Omega)_{1/2}} + \|u(t_2)\|_{DL_p(\Omega)_{1/2}})\|u(t_1)-u(t_2)\|_{DL_p(\Omega)_{1/2}}$$
$$\leq ML_0^2 C(\|u(t_1)\|_{DL_p(\Omega)_{1/2}} + \|u(t_2)\|_{DL_p(\Omega)_{1/2}})|t_1-t_2|^\beta. \tag{12}$$

We used the Lemma 6 in the above third step and the formula (10) in fifth step.

( Step 3) Suppose that $u_0 \in DL_p(\Omega)_{1/2}$ is smooth satisfying $\dfrac{\partial u_{0i}}{\partial x_j}=0$ $(i \neq j)$. Then $u_0$ is divergence free from lemma 1 and $(t_0, u_0)\in U$. Set

$$V = B_\varepsilon(t_0,u_0) = \left\{(t,u(t))\in U : |t-t_0|\langle \varepsilon \leq 1, \|u-u_0\|_{DL_p(\Omega)_{1/2}}\langle \varepsilon\right\}.$$

Then for $(t,u(t))\in V$,

$$\|u\|_{DL_p(\Omega)_{1/2}} = \|u-u_0+u_0\|_{DL_p(\Omega)_{1/2}} \leq \|u-u_0\|_{DL_p(\Omega)_{1/2}} + \|u_0\|_{DL_p(\Omega)_{1/2}} \leq \varepsilon + \|u_0\|_{DL_p(\Omega)_{1/2}}.$$

Let $L = \varepsilon + \|u_0\|_{DL_p(R^3)_{1/2}}$, $L_1 = 2ML_0^2 L$, $L_2 = 2ML_0^2 L(C_1+C)$, $L_3 = Max(L_2,L_1)$ and $\beta_2 = Min(\beta,\beta_1)$, then from (10),(11) and (12) for all $(t_i,u_i)\in V$ we have

$$\|F(t_1,u_1(t_1))-F(t_2,u_2(t_2))\|_{DL_p(\Omega)}$$
$$\leq \|F(t_1,u_1(t_1))-F(t_1,u_2(t_1))\|_{DL_p(\Omega)} + \|F(t_1,u_2(t_1))-F(t_2,u_2(t_2))\|_{DL_p(\Omega)}$$
$$= \|(u_1(t_1)\bullet\nabla)u_1(t_1)-(u_2(t_1)\bullet\nabla)u_2(t_1)\|_{DL_p(\Omega)} + \|(u_2(t_1)\bullet\nabla)u_2(t_1)-(u_2(t_2)\bullet\nabla)u_2(t_2)\|_{DL_p(\Omega)}$$
$$\leq 2ML_0^2 L\|u_1(t_1)-u_2(t_2)\|_{DL_p(\Omega)_{1/2}} + 2ML_0^2 LC_1|t_1-t_2|^{\beta_1} + 2ML_0^2 LC|t_1-t_2|^\beta$$
$$\leq L_1\|u_1(t_1)-u_2(t_2)\|_{DL_p(\Omega)_{1/2}} + L_2|t_1-t_2|^{\beta_2}$$
$$\leq L_3(|t_1-t_2|^{\beta_2} + \|u_1(t_1)-u_2(t_2)\|_{DL_p(\Omega)_{1/2}}).$$

Hence $F(t,u(t))$ satisfies the assumption $(F)$. Therefore from lemma 8 for every initial data

$(t_0, u_0) \in U_1$ the initial value problem (9) has unique local solution

$$u(t) \in C([t_0, t_1) : DL_p(\Omega)) \cap C^1((t_0, t_1) : DL_p(\Omega)) \tag{13}$$

where $t_1 = t_1(u_0)$.

The solution (13) of (9) is also the solution of (4). The Theorem 3.4 in [5] mean that as long as the solution of (4) exists, this solution is smooth. From Theorem 3.4 in [5] we have the solution $u(t,x) \in \left(C^\infty\left([0,t_1) \times \overline{\Omega}\right)\right)^3$. Substituting $u(t,x)$ into (1.1) we get the solution $p(t,x)$. We also have $p(t,x) \in C^\infty\left([0,t_1) \times \overline{\Omega}\right)$ It follows from the lemma 1 and $u \in DL_p(\Omega)$ that the solution $u(t,x)$ is divergence-free. Changing the value of $u$ on $\partial\Omega$ to zero we get a unique smooth local strong solution of the Navier-Stokes initial value problem (1.1)-(1.4). □

If $\Omega$ is a bounded open subset of $R^3$ then the solution (13) is a classical solution of (9) and (3.1)-(3.2).

### Acknowlements

The first author is grateful to her supervisors Chi-Kun Lin and Xinyao Yang for their teaching and cultivation.

### References

[1] S.Agmon, On the eigenfunctions and on the eigenvalues of general elliptic value problems. Comm.Pure Appl.Math. 15(1962)119-147.

[2] S.Agmon and L.Nireberg, Properties of solutions of ordinary differential equations in Banach spaces, Comm.Pure Appl.Math. 16(1963)121-239.

[3] Andras Batkai, Marjeta Kramar Fijavz, Abdelaziz Rhandi, Positive Operator Semigroups, Birkhauser (2017).

[4] H.Fujita,T.Kato, On the Navier-Stokes initial value problem I, Arch. Ration. Mech. Anal. 16(1964),269-315.

[5] Y.Giga,T.Miyakava, Solutions in $L_r$ of the Vavier-Stokesinitial value problem,

Arch.Ration.Mech.Anal. 89(1985)267-281.

[6] Philip Isett，Hölder Continuous Euler Flows in Three Dimensions with Compact Support in Time，Princeton University Press,(2017).

[7] T.Kato, Strong $L_p$-solution of the Navier-Stokes equation in $R^m$, with application to weak solutions, Math. Z. 187(1984)471-480.

[8] Engel Klaus-Jochen, Rainer Nagel. One-Parameter Semigroups for Linear Evolution Equations, Springer, (2000).

[9] E.Kreyszig, Introduction functional analysis with applications, John Wiley &Sons (1978).

[10] Peter Meyer-Nieberg, Banach Lattices, Springer-Verlag, (1991).

[11] Mukhtarbay.Otelbaev, Existence of a strong solution of the Navier-Stokes equation, Mathematical Journal( Russia), 13(4)(2013),5-104.

[12] A.Pazy, Semigroups of linear operators and applications to partial differential equations, Springer Verlag (1983,reprint in China in 2006).

[13] James C. Robinson, Infinite-Dimentional Dynamical Systems, Cambridge University Press, (2001).

[14] Walter Rudin, Principals of Mathematical Analysis, China Machine Press, 2004.

[15] Walter Rudin, Real and Complex analysis, China Machine Press, 2004

[16] Veli B.Shakhmurov, Nonlocal Navier-Stokes problems in abstract function, Nonlinear Analysis:Real World Applications , 26(2015)19-43.

[17] Maoting Tong, Daorong Ton, A local strong solution of the Navier-Stokes problem in $L_2(\Omega)$, Journal of Mathematical Sciences: Advances and applications,

Volume No:62, Issue no:1,May(2020).1-19. Available at http://scientificadvances.co.in

DOI: http://dx.doi.org/10.18642/jmsaa_7100122125